\documentclass[12pt]{article}
\usepackage{amsxtra,amssymb,amsmath,enumerate,graphicx,bbm,hyperref,amsthm,color}
\usepackage[active]{srcltx}
\usepackage{t1enc}
\usepackage{aurical}
\usepackage[T1]{fontenc}
\usepackage{float}
\usepackage[linesnumbered,ruled,vlined]{algorithm2e}
\usepackage{subcaption}
\usepackage{titlesec}
\usepackage{float}

\newtheorem{lemma}{Lemma}[section]

\newtheorem{theorem}[lemma]{Theorem}

\newcommand{\e}{\mathbb{E}}
\newcommand{\E}{\mathbb{E}}

\renewcommand{\P}{\mathbb{P}}
\newcommand{\R}{\mathbb{R}}

\newcommand{\1}{{\mathbf 1}}

\def\be{\begin{eqnarray}}
\def\ee{\end{eqnarray}}
\def\b*{\begin{eqnarray*}}
\def\e*{\end{eqnarray*}}

\def\vp{\varphi}

\def\be{\begin{eqnarray}}
\def\ee{\end{eqnarray}}
\def\beq{\begin{equation}}
\def\eeq{\end{equation}}
\def\b*{\begin{eqnarray*}}
\def\e*{\end{eqnarray*}}
\def\bi{\begin{itemize}}
\def\ei{\end{itemize}}

\def \1{{\bf 1}}
\def\vp{\varphi}

\def\={\;=\;}

\def\x{\times}

\def\esssup{{\rm ess}\!\sup\limits}

\def \proof{{\noindent \bf Proof. }}
\def \ep{\hbox{ }\hfill$\Box$}
 \def\reff#1{{\rm(\ref{#1})}}

 \def\vs#1{\vspace{#1mm}}

\def\ti{{t_i}}
\def\tip{ {t_{i+1}} }

\def \E{\mathbb{E}}
\def \F{\mathbb{F}}

\def \N{\mathbb{N}}
\def \P{\mathbb{P}}

\def \R{\mathbb{R}}
\def\T{\mathbb{T}}

\def\Bc{{\cal B}}

\def\Dc{{\cal D}}

\def\Fc{{\cal F}}

\def\Kc{{\cal K}}

\def\Pc{{\cal P}}

\def\Wc{{\cal W}}

\def\Lb{{\mathbf L}}

\def\Yb{\bar Y}
\def\Zb{\bar Z}

\def\Yt{\tilde Y}
\def\Zt{\tilde Z}

\newtheorem{Proposition}{Proposition}[part]

\newtheorem{Lemma}{Lemma}[part]

\newtheorem{Remark}{Remark}[part]

\makeatletter \@addtoreset{equation}{section}

\@addtoreset{Definition}{section}

\@addtoreset{Theorem}{section}

\@addtoreset{Proposition}{section}

\@addtoreset{Assumption}{section}

\@addtoreset{Corollary}{section}

\@addtoreset{Lemma}{section}

\@addtoreset{Remark}{section}

\@addtoreset{Example}{section}

\def\Sb{\mathbf{S}}

\def\Xbf{\mathbf{X}}

\def\Rc{{\cal R}}
\def\ur{{\rm u}}

\begin{document}

\title{Numerical approximation of general Lipschitz BSDEs with branching processes}

\author{Bruno Bouchard\footnote{Universit\'e Paris-Dauphine, PSL Research University, CNRS, UMR [7534], CEREMADE, 75016 PARIS, FRANCE.}  \thanks{bouchard@ceremade.dauphine.fr} 
\and Xiaolu Tan\addtocounter{footnote}{-1}\footnotemark[\value{footnote}]~\addtocounter{footnote}{1}\thanks{tan@ceremade.dauphine.fr}
\and Xavier Warin\footnote{EDF R\&D \& FiME, Laboratoire de Finance des March\'es de l'Energie} \thanks{warin@edf.fr} 
}

\maketitle

\begin{abstract}
	We extend the branching process based numerical algorithm of Bouchard et al.~\cite{BouchardTanZou}, that is dedicated to semilinear PDEs (or BSDEs) with Lipschitz nonlinearity, to the case where the nonlinearity involves the gradient of the solution. As in  \cite{BouchardTanZou}, this requires a localization procedure that uses a priori estimates on the true solution, so as to ensure the well-posedness of the involved Picard iteration scheme, and the global convergence of the algorithm. When, the nonlinearity depends on the gradient, the later needs to be controlled as well. This is done by using a face-lifting procedure.   
	Convergence of our algorithm is proved  without any limitation on the time horizon.
	We also provide numerical simulations  to illustrate the performance of the algorithm.
\end{abstract}

\vspace{2mm}

\section{Introduction}

	The aim of this paper is to extend the branching process based numerical algorithm proposed in Bouchard et al.~\cite{BouchardTanZou} to general BSDEs in form:
	\be\label{eq:BSDE}
		Y_{\cdot}=g({X}_T)+\int_{\cdot}^{T}f( {X}_s,Y_s,  Z_s)\,ds-\int_{\cdot}^{T}  Z_s^{\top} dW_s,
	\ee
	where $W$ is a standard $d$-dimensional Brownian motion, $f: \R^d \x \R \x \R^d \to \R$ is the driver function, 
	$g: \R^d \to \R$ is the terminal condition, and $X$ is the solution of 
	\begin{equation} \label{eq: Diffusion}
		{X}=X_{0}+\int_{0}^{\cdot} \mu({X}_s)\,ds+\int_{0}^{\cdot} \sigma( {X}_s)\,dW_s,
	\end{equation}
	with constant initial condition $X_0 \in \R^{d}$ and coefficients  $(\mu,\sigma):\R^{d}\mapsto \R^{d}\x \R^{d\x d}$, that are assumed to be Lipschitz\footnote{As usual, we could add a time dependency in the coefficients $f$, $\mu$ and $\sigma$ without any additional difficulty. }.
	From the PDE point of view, this amounts to solving the parabolic equation
	\begin{equation*} \label{eq:PDEsemiLinear}
		 \partial_t u + \mu \cdot Du + \frac{1}{2} {\rm Tr}[\sigma \sigma^{\top}  D^2 u]
		+ f\big(\cdot, u, Du\sigma \big)     = 0,
		~~~u(T,\cdot) = g.
	\end{equation*}
		
	The main idea of \cite{BouchardTanZou} was to  approximate the driver function by   local polynomials  and use a Picard iteration argument so as to reduce the problem to solving BSDE's with (stochastic) global polynomial drivers, see Section \ref{sec:model}, to which the branching process based pure forward Monte-Carlo algorithm of  \cite{Henry-Labordere_branching, HOTTW, HTT} can be applied. See for instance \cite{McKean_1975,Skorokhod,Watanabe} for the related Feynman-Kac representation of the KPP (Kolmogorov-Petrovskii-Piskunov) equation. 
	
	This  algorithm   seems to be very adapted to situations where the original driver can be well approximated by polynomials with rather small coefficients on quite large domains. The reason is that, in such a situation, it is basically a pure forward Monte-Carlo method, see  in particular \cite[Remark 2.10(ii)]{BouchardTanZou}, which can be expected to be less costly  than the classical schemes, see e.g.~\cite{BallyPages,BouchardTouzi,BouchardWarin,GLW,Zhang} and the references therein. However, the numerical scheme  of \cite{BouchardTanZou} only works when the driver function $(x,y,z)\mapsto f(x,y,z)$ is independent of $z$, i.e.~the nonlinearity in the above equation does not depend on the gradient of the solution. 

	Importantly, the algorithm proposed in \cite{BouchardTanZou} requires the truncation of the approximation of the $Y$-component at some given time steps. The reason is that BSDEs with polynomial drivers may only be defined up to an explosion time. This truncation is based on a priori estimates of the true solution. It ensures  the well-posedness of the algorithm on an arbitrary time horizon, its stability, and global convergence. 
	
	In the case where the driver also depends on the $Z$ component of the BSDE, a similar truncation has to be performed on the gradient itself. It can however not be done by simply projecting $Z$ on a suitable compact set at certain time steps, since $Z$ only maters up to an equivalent class of $([0,T]\x \Omega,dt\x d\P)$. Alternatively, we propose to use a face-lift procedure at certain time steps, see 
	\eqref{eq: facelift}. Again this time steps depend on the explosion times of the corresponding BSDEs with polynomial drivers. Note that a similar face-lift procedure is used in Chassagneux, Elie and Kharroubi \cite{contrainteZ}\footnote{We are grateful to the authors for several discussions on this subject.} in the context of the discrete time approximation of BSDEs with contraint on the $Z$-component. 
\\
	
	We prove the convergence of the scheme.  The very good performance of this approach is illustrated in  Section \ref{sec:NumericalExamples} by a numerical test case. 
\\

{\bf Notations: } All over this paper, we view elements of $\R^{d}$, $d\ge 1$, as column vectors. Transposition is denoted by the superscript ${}^{\top}$. We consider a complete filtered probability space  $(\Omega,\Fc,\F=(\Fc_{t})_{t\le T},\P)$ supporting a $d$-dimensional Brownian motion $W$. We simply write $\E_{t}[\cdot]$ for $\E[\cdot|\Fc_{t}]$, $t\le T$.   We use the standard notations $\Sb_{2}$ (resp.~$\Lb_{2}$) for the class of progressively measurable processes $\xi$ such that $\|\xi\|_{\Sb_{2}}:=\E[\sup_{[0,T]}\|\xi\|^{2}]^{\frac12}$ (resp.~$\|\xi\|_{\Lb_{2}}:=\E[\int_{0}^{T}\|\xi_{s}\|^{2}ds]^{\frac12}$) is finite. The dimension of the process $\xi$ is given by the context. For a map $(t,x)\mapsto \psi(t,x)$, we denote by $\partial_{t} \psi$ is derivative with respect to its first variable and by $D\psi$ and $D^{2} \psi$ its Jacobian and Hessian matrix with respect to its second component.

\section{Approximation of BSDE using local polynomial drivers and the Picard iteration}
\label{sec:model}

	For the rest of this paper, let us consider the (decoupled) forward-backward system \eqref{eq:BSDE}-\eqref{eq: Diffusion} in which  $f$ and $g$ are both bounded  and Lipschitz continuous, and 
	$\sigma$ is non-degenerate such that there is a constant $a_{0}>0$ satisfying
	\begin{align}\label{eq: unif elliptic}
		\sigma \sigma^{\top}(x) \ge |a_0|^{2} {\rm I}_d, ~\forall x \in \R^d.
	\end{align}
	We also assume that  $\mu$, $\sigma$,   $D \mu$ and $D \sigma$  are  all bounded and continuous.  
	In particular,  \eqref{eq:BSDE}-\eqref{eq: Diffusion} has a unique solution $(X,Y,Z)\in \Sb_{2}\x \Sb_{2}\x \Lb_{2}$.   The above conditions indeed imply that $|Y| + \|\sigma(X)^{-1}Z\| \le M \mbox{ on } [0,T]$, for some $M>0$.
 
  \begin{Remark} The above assumptions can be relaxed by using standard localization or mollification arguments. For instance, one could simply assume that $g$ has polynomial growth and is locally Lipschitz. In this case, it can be truncated outside a compacted set so as to reduce to the above.  Then, standard estimates and stability results for SDEs and BSDEs can be used to estimate the additional error in a very standard way.  See e.g.~\cite{ElKarouiBSDE}.
 \end{Remark}

\subsection{Local polynomial approximation of the generator}\label{sec: local poly approx}

	As in \cite{BouchardTanZou}, our first step is to approximate the driver $f$ by a driver $f_{{\circ}}$ 
	that has a local polynomial structure. The difference is that it now depends on both components of the solution of the BSDE. 
	Namely, let 
	\be\label{eq: f loc poly} 
		f_{\circ}(x,y,z, y', z'):= \sum_{j=1}^{j_{\circ}} \sum_{\ell \in L} c_{j,\ell}(x)y^{\ell_{0}} \prod_{q=1}^{q_{\circ}} (b_q(x)^{\top} z)^{\ell_q} \vp_{j}(y', z'),
	\ee
	in which $(x,y,z, y', z')\in   \R \x \R\x \R^d \x \R \x \R^d$,
	$L:= \{0,\cdots,L_{\circ}\}^{q_{\circ}+1}$ (with $L_{\circ},q_{\circ}\in \N$),
	the functions $(b_q)_{0\le q \le q_{\circ}}$ and $(c_{j,\ell},\vp_{j})_{\ell \in L,1\le  j\le j _{\circ}}$ (with $j_{\circ}\in \N$) are   continuous and satisfy 
	\be\label{eq: borne a vp}
		|c_{j,\ell}  | \le C_L \;,\; 
		\|b_q\| \le 1,\;\;
		|\vp_{j} |\le 1,
		\mbox{ and } \|D \vp_{j}\| \le L_{\varphi}, 
	\ee 
	for all $1\le j\le j_{\circ}$, $0\le q\le q_{\circ}$ and $\ell \in L$,  for some constants $C_{L},L_{\varphi}\ge 0$.

	For a good choice of the local polynomial $f_{\circ}$, we can assume that	 
	\b* \label{eq:f_0}
		 (x, y,z) \mapsto 		 
		 \bar f_{\circ} \big(x, y,  z\big):=f_{\circ} \big(x, y,  z, y, z\big)
	\e*
	is globally bounded and Lipschitz.
	Then, the BSDE
\be\label{eq:BSDElo}
	\Yb_{\cdot}
	=
	g({X}_T)
	+\!
	\int_{\cdot}^{T} \bar f_{{\circ}}( {X}_s,\Yb_s, \Zb_s)\,ds
	\!-\!\int_{\cdot}^{T} \Zb_s^{\top}\,dW_s,
\ee
has a unique solution $( \Yb,  \Zb)\in \Sb_{2}\x \Lb_{2}$, and  standard estimates imply that $(\Yb, \Zb)$ provides a good approximation of $(Y,Z)$ whenever $\bar f_{{\circ}}$ is a good approximation of $f$:
\be \label{eq: erreur driver}
	\|Y- \Yb \|_{\Sb_{2}}+\|Z- \Zb\|_{\Lb_{2}}\le
	C_{\circ} ~\|(f-\bar f_{{\circ}})(X,Y ,Z)\|_{\Lb_{2}},
\ee
for some $C_{\circ}>0$ that depends on the global Lipschitz constant of $\bar f_{\circ}$ (but not on the precise expression of $\bar f_{{\circ}}$), see e.g.~\cite{ElKarouiBSDE}.

One can think at the $(c_{j,\ell})_{\ell \in L}$ as the coefficients of a polynomial 
approximation of $f$ in terms of $(y, (b_q(x)^{\top} z))_{q\le q_{\circ}}$ on a subset $A_{j}$, the $A_{j}$'s forming a partition of $[-M,M]^{1+d}$. Then, the $\vp_{j}$'s have to be considered as smoothing kernels that allow one to pass in a Lipschitz way from one part of the partition to another one.

The choice of the basis functions $(b_q)_{q \le q_{\circ}}$ as well as $(\vp_{j})_{1\le j\le j_{\circ}}$ will obviously depend on the application, but it should in practice typically be constructed  
such that the sets 
\be\label{eq: def Sj}
A_{j}:=\{(y,z) \in \R\x\R^{d} : \vp_{j}(y,z) = 1\}
\ee
 are large and the intersection between the supports of the $\vp_{j}$'s are small.  See \cite[Remark 2.10(ii)]{BouchardTanZou} and below.
Finally, since the function $\bar f_{\circ}$ is chosen to be globally bounded and Lipschitz, by possibly adjusting the constant $M$, we can assume without loss of generality that
\begin{align}\label{eq: bone bar Y}
	|\bar Y| + \|\bar Z^{\top} \sigma(X)^{-1}\| ~\le~ M.
\end{align}

For later use, let us recall that $\bar Y$ is related to the unique bounded and continuous viscosity solution  $\bar u$ of 
	\begin{equation*} \label{eq:PDEsemiLinear}
		 \partial_t \bar u + \mu\cdot  D\bar u + \frac{1}{2} {\rm Tr}[\sigma \sigma^{\top}  D^2 \bar u]
		+ \bar f_{\circ}\big(\cdot, \bar u, D\bar u\sigma \big)     = 0,
		~~~\bar u(T,\cdot) = g,
	\end{equation*}
	though 
	\be\label{eq: Yb:bar u}
	\bar Y=\bar u(\cdot,X). 
	\ee
Moreover, 
\b*\label{eq: bar u borne lip}
\mbox{$\bar u$ is bounded by $M$ and $M$-Lipschitz.}
\e*
	
\subsection{Picard iteration with truncation and face-lifting}

	Our next step is to introduce a Picard iteration scheme to approximate the solution $\Yb$ of  \eqref{eq:BSDElo} so as to be able to apply the branching process based forward Monte-Carlo approach of \cite{Henry-Labordere_branching, HOTTW, HTT} to each iteration: given $(\Yb^{m-1},\Zb^{m-1})$, use the representation of    the BSDE with driver $f_{\circ}(X,\cdot,\cdot,\Yb^{m-1},\Zb^{m-1})$.

		However, although  the map $(y,z) \mapsto \bar f_{\circ}(x,y,z)=f_{\circ}(x, y,z, y,z)$ is globally Lipschitz,
	the map $(y,z) \mapsto f_{\circ}(x, y,z, y',z')$ is a polynomial, given fixed $(x,y', z')$,
	and hence only locally Lipschitz in general.
	In order to reduce to a Lipschitz driver, we need to truncate the solution at certain time steps, that are smaller than the natural explosion time of the corresponding BSDE with (stochastic) polynomial driver. 
	As in  \cite{BouchardTanZou}, it can be performed by a simple truncation at the level of the first component of the solution. As for the second component, that should be interpreted as a gradient, a simple truncation does not make sense, the gradient needs to be modified by modifying the function itself. Moreover, from the BSDE viewpoint, $Z$ is only defined up to an equivalent class on $([0,T]\x\Omega,dt\x d\P)$, so that changing its value at a finite number of given times does not change the corresponding BSDE. We instead use a face-lifting procedure,   as in \cite{contrainteZ}. 
	\vs2
	
	More precisely, let us define the operator $\Rc$   on the space of bounded functions $\psi :\R^{d}\mapsto \R$ by 
	$$
	\Rc[\psi]:=(-M)\vee \left[\sup_{p \in \R^{d}} \big( \psi(\cdot+p) - \delta(p) \big)\right]\wedge M,
	$$ 
	where 
		 $$
		 \delta(p) := \sup_{q \in [-M, M]^d} (p^{\top} q) = M \sum_{i=1}^d |p_i|.
		 $$
 	The operation $\psi\mapsto \sup_{p \in \R^{d}} \big( \psi(\cdot+p) - \delta(p) \big)$ maps $\psi$ into the smallest $M$-Lipschitz function above $\psi$. This is the so-called face-lifting procedure, which has to be understood as a (form of) projection on the family of $M$-Lipschitz functions, see e.g.~\cite[Exercise 5.2]{BouchardChassagneux}, see also Remark \ref{rem: facelift non exp} below. 	The outer operations in the definition of $\Rc$ are just a natural projection on $[-M,M]$.  
\vs2

	Let now $(h_{\circ},M_{h_{\circ}})$ be such that \eqref{eq: def Mho et ho} and \eqref{eq: def Mho et ho Z}  in the Appendix hold. The constant $h_{\circ}$ is  a lower bound for the explosion time of the BSDE with driver $(y,z) \mapsto f_{\circ}(x, y, z, y', z')$ for any fixed $(x,y',z')$. Let us then fix $h \in (0, h_{\circ})$ such that $N_{h}:=T/h \in \N$, and define  
	\begin{align*} 
		\ti=ih
		\;\;\;\;\;\mbox{and}\;\;\;\;
		\T := \{\ti, ~i= 0, \cdots, N_{h}\}.
	\end{align*}
	Our algorithm consists in using a Picard iteration scheme to solve \eqref{eq:BSDElo}, which re-localize the solution at each time step of the grid $\T$ by applying operator $\Rc$.
	
	Namely, using the notation  $X^{t,x}$ to denote  the solution of \eqref{eq: Diffusion} on $[t,T]$ such that $X^{t,x}_{t}=x\in \R^{d}$, we initialize our  Picard scheme by setting  
	\begin{align*}
		Y^{x,0}_{T}=\Yb^{x,0}_{T}&=g(x)\\
		(Y^{x,0},Z^{x,0})&=(\Yb^{x,0},\Zb^{x,0}) =({\rm y},D{\rm y})(\cdot, X^{t_{i},x})\;\mbox{ on }  [t_{i},t_{i+1}),\;i\le N_{h}-1, 
	\end{align*}
	in which ${\rm y}$ is a continuous function, $M$-Lipschitz in space,  continuously differentiable in space on $[0,T)\x \R^{d}$ and such that  ${\rm |y|} \le M$ and ${\rm y}(T,\cdot)=g$.
	Then, given  $(\Yb^{x,m-1},\Zb^{x,m-1})$, for $m\ge 1$, we define $(\Yb^{x,m}, \Zb^{x,m})$ as follows:
	\begin{enumerate}
		\item For $i = N_{h}$, set $\bar \ur^{m}_{t_{i}}=\bar \ur^{m}_{T}:=g$ 
		\item For $i<N_{h}$, given $( \Yb^{x,m-1}, \Zb^{x,m-1})$:
		\begin{enumerate}
		\item Let $(Y^{x,m}_{\cdot}, Z^{x,m}_{\cdot})$ be the unique solution on $[t_{i},t_{i+1})$ of  
		\begin{align}  
			Y^{x,m}_{\cdot} 
			&=
			\bar \ur^{m}_{t_{i+1}}(X^{t_{i},x}_{t_{i+1}}) \nonumber
			 \\
			 &+
			 \int^{t_{i+1}}_{\cdot} \!\!\! f_{\circ} 
			\big(X^{t_{i},x}_s, Y^{x,m}_s, Z^{x,m}_s, \Yb^{x,m-1}_s,   \Zb^{x,m-1}_s \big) ds 
			\nonumber\\
			&- \int^{t_{i+1}}_{\cdot} \!\!\!\! (Z^{x,m}_s)^{\top} dW_s.
			\label{eq:def_Ym}
		\end{align}

		\item \label{item:face_lifting} 
		Let $\ur^{m}_{t_{i}} :   x\in \R^d \mapsto Y^{x,m}_{t_{i}}$,
		and set 
		\be\label{eq: facelift}
		\bar \ur^{m}_{t_i} := \Rc[\ur^{m}_{t_i}].
		\ee
		\item Set  $\Yb^{x,m} := Y^{x,m}$ on $ (t_i, t_{i+1}]$,  $\Yb^{x,m}_{t_i} :=  \bar \ur^{m}_{t_i}(x)$,
		and $ \Zb^{x,m} := Z^{x,m}$ on $[t_i, t_{i+1})$, for $x\in \R^{d}$.
		\end{enumerate}
		\item We finally define  $\Yb^{m}_{T}=g(X_{T})$ and
		\be \label{eq:def_Yb_m}
			 (\Yb^{m},\Zb^{m})
			:=
			(\Yb^{X_{t_{i}},m},\Zb^{X_{t_{i}},m})\;\mbox{ on } [t_{i},t_{i+1}), \;i\le N_{h}.
		\ee
	\end{enumerate}
	
	In above, the existence and uniqueness of the solution $(Y^{x,m}, Z^{x,m})$ to \eqref{eq:def_Ym} is ensured by Proposition \ref{prop: bound tilde Y tilde Z}.
	The projection operation in \eqref{eq: facelift} is consistent with the behavior of the solution of  \eqref{eq:BSDElo}, recall \eqref{eq: bone bar Y}, and it is crucial to control the explosion of $(\Yb^{m},\Zb^{m})$ and therefore to ensure both the stability and the convergence of the scheme. This procedure is non-expansive, as explained in the following Remark, and therefore can not alter the convergence of the scheme. 
\begin{Remark}\label{rem: facelift non exp} Let $\psi,\psi'$ be two measurable and bounded maps on $\R^{d}$.  Then, $ \sup_{p \in \R^{d}} |\psi(\cdot+p)-\delta(p)-\psi'(\cdot+p)+\delta(p)|$ $=$ $\sup_{x \in \R^{d}} |\psi(x) -\psi'(x)|$, and therefore 
$\|\Rc[ \psi]-\Rc[\psi']\|_{\infty}$ $\le$ $ \|  \psi-  \psi'\|_{\infty}$. In particular, since  $\bar u$ defined through \eqref{eq: Yb:bar u} is $M$-Lipschitz in its space variable and bounded by $M$, we have  $\bar u(t,\cdot) =\Rc[\bar u(t,\cdot)]$ for $t\le T$ and therefore 
$$
\|\Rc[\psi]-\bar u(t,\cdot)\|_{\infty}\le  \|  \psi-  \bar u(t,\cdot)\|_{\infty}
$$
for all $t\le T$ and all measurable and bounded map  $\psi$.
\end{Remark}
 	
Also note that, if we had  $(\Yb^{m-1}_{t},\Zb^{m-1}_{t})\in A_{j}$ if and only if  $(\Yb_{t},\Zb_{t})\in A_{j}$, for all $j\le j_{\circ}$, then we would have $(\Yb^{m-1},\Zb^{m-1})=(\Yb,\Zb)$, recall \eqref{eq: def Sj} and the definition of $\bar f_{\circ}$ in terms of $f_{\circ}$. This means that we do not need to be very precise on the original prior, whenever the sets $A_{j}$ can be chosen to be large.

	\vs2
	
	From the theoretical viewpoint, the error due to the above Picard iteration scheme can be deduced from classical arguments. Recall that $(h_{\circ},M_{h_{\circ}})$ is such that \eqref{eq: def Mho et ho} and \eqref{eq: def Mho et ho Z} in the Appendix hold. 
	
	\begin{theorem}\label{thm: main} 
		For each $m\ge 0$, the algorithm defined in {\rm 1.-2.-3.~}above provides  the unique solution $(\Yb^{m}, \Zb^m)\in \Sb_{2}\x \Lb_{2}$.
		Moreover, it satisfies     
		$|\bar Y^m|\vee \|\bar Z^{m\top}\sigma(X)^{-1}\|  \le M_{h_{\circ}}$, and there exists a measurable map $(\bar u^{m},\bar v^{m}): [0,T]\x \R^{d}\mapsto \R^{1+d}$, that is continuous on  $\cup_{i<N_{h}}(t_{i},t_{i+1})\x \R^{d}$,  such that $\bar u^{m}(t_{i},\cdot)$ is continuous on $\R^{d}$ for all $i\le N_{h}$, and 
		\begin{align}\label{eq: def bar um bar vm}
		\bar Y^{m}&=\bar u^{m}(\cdot, X) \mbox{ on $[0,T]$ $\P$-a.s.}\\
		\bar Z^{m}&=\bar v^{m}(\cdot, X)\mbox{ $dt\x d\P$-a.e.~on $[0,T]\x \Omega$.}\nonumber\\
		\bar v^{m \top}&=D\bar u^{m}\sigma  \mbox{ on $\cup_{i< N_{h}}(t_{i},t_{i+1})\x \R^{d}$.}\nonumber
		\end{align}		
		Moreover, for any constant $\rho \in (0,1)$, there is some constant $C_{\rho} > 0$ such that
		\b*
			 |\Yb^{m}_{t}- \Yb_{t} |^2 +\E_{t}[\int_{t}^{T}\|\Zb^{m}_{s}-\Zb_{s}\|^{2}ds]
			~\le~
			C_{\rho}\; \rho^m,\;\mbox{ for all $t\le T$.}
		\e*
	\end{theorem}
	\proof
	$\mathrm{i)}$ Recall from Remark \ref{rem: facelift non exp} that $\Rc$ maps bounded functions into $M$-Lipschitz functions that are bounded by $M$.
	Then, by Proposition \ref{prop: bound tilde Y tilde Z} in the Appendix,
	the solutions $(Y^{x,m}, Z^{x,m})$ as well as $(\bar Y^{x,m}, \bar Z^{x,m})$ are uniquely defined in and below \eqref{eq:def_Ym}.
	Moreover, one has $|\bar Y^{x,m}_{t}|\vee \|(\bar Z^{x,m}_{t})^{\top}\sigma(X^{t_{i},x}_{t})^{-1}\| $ $\le$ $M_{h_{\circ}}$,  for all $x\in \R^{d}$, $i<N_{h}$ and $t\in [t_{i},t_{i+1})$.
	As a consequence,   $(\bar Y^{m},\bar Z^{m})$ is also uniquely defined  and satisfies $| \bar  Y^m| \vee \| \bar Z^{m\top}\sigma(X)^{-1}\| \le  M_{h_{\circ}}$.  
	Using again Proposition \ref{prop: bound tilde Y tilde Z}, one has the existence of $(\bar u^m, \bar v^m)$ satisfying the condition in the statement.

	\vspace{1mm}

	$\mathrm{ii)}$  We next prove the convergence of the sequence $(\bar Y^m, \bar Z^m)_{m \ge 0}$ to $(\bar Y, \bar Z)$.
	Since $\{(\bar Y^{x,m}, \bar Z^{x,m}), x\in \R^{d}\}$ is uniformly bounded, 
	the generator $f_{{\circ}}(x,y,$ $ z,$ $y',$ $z')$ in \eqref{eq:def_Ym} can be considered to be uniformly Lipschitz in $(y,z)$ and $(y', z')$.
	Assume that the corresponding Lipschitz constants are $L_1$ and $L_2$.	
	
	\vspace{1mm}

	Let us set  $\Theta^{x,m}:=(Y^{x,m},Z^{x,m})$ and define $(\Delta Y^{x,m},\Delta Z^{x,m}) := (Y^{x,m} - \Yb^{x},Z^{x,m} - \Zb^{x})$ where $\bar \Theta^{x}:=(\Yb^{x},\Zb^{x})$ denotes the solution of 
	\begin{align*} 
			\Yb^{x}_{\cdot} 
			&=
			\bar u(t_{i+1},X^{t_{i},x}_{t_{i+1}})
			 +
			 \int^{t_{i+1}}_{\cdot} \!\!\! \bar f_{\circ} 
			\big(X^{t_{i},x}_s, \Yb^{x}_s, \Zb^{x}_s \big) ds 
			- \int^{t_{i+1}}_{\cdot} \!\!\!\! (\Zb^{x}_s)^{\top} dW_s,\nonumber
		\end{align*}
		on each $[t_{i},t_{i+1}]$, recall \eqref{eq: Yb:bar u}.
	In the following, we fix $\beta>0$ and use the notation
	$$
	\|\xi\|_{\beta, t}:=\E_{t}[\int_{t}^{t_{i+1}}e^{\beta s}|\xi_{s}|^{2}ds]^{\frac12}\;\mbox{ for } \xi \in \Lb_{2},\;t\in [t_{i},t_{i+1}).
	$$
  Fix $t\in [t_{i},t_{i+1})$. By applying It\^o's formula to $(e^{\beta s} (\Delta Y^{x,m+1}_{s})^2)_{s\in [t,t_{i+1}]}$ 
	and then taking expectation, we obtain
	\begin{align*}
		&
		e^{\beta t}   |\Delta Y^{x,m+1}_{t}|^2
		+   \beta \| \Delta Y^{x,m+1}\|_{\beta, t}^2 
		+ 2 \| \Delta Z^{x,m+1}\|_{\beta, t}^2  \\
		&\le
		\E_{t} \big[ e^{\beta t_{i+1}} (\Delta Y^{x,m+1}_{t_{i+1}-})^2 \big]
		\\
		&+
		2 \E_{t} \Big[ \int_{t}^{t_{i+1}}e^{\beta s} \Delta Y^{x,m+1}_s \big( f_{{\circ}} (X^{t_{i},x}_s,   \Theta^{x,m+1}_s,  \Theta^{x,m}_s) 
		 - f_{{\circ}}(X^{t_{i},x}_s,\bar \Theta^{x}_s, \bar \Theta^{x}_s) \big) ds \Big].
	\end{align*}
	Using the Lipschitz property of $f_{\circ}$ and the inequality $\lambda + \frac{1}{\lambda} \ge 2$ for all $\lambda >0$,
	it follows that, for all $\lambda_1, \lambda_2 > 0$,
	\begin{align} \label{eq:contract_beta}
		&e^{\beta t}   |\Delta Y^{x,m+1}_{t}|^2
		+ (\beta - (2L_1 + \lambda_1 L_1 +  \lambda_2 L_2)) \| \Delta Y^{x,m+1} \|_{\beta, t}^2
		\\
		&+ (2 - \frac{L_1}{\lambda_1}) \| \Delta Z^{x,m+1} \|_{\beta, t}^2 \nonumber \\
		&\le
		\E_{t} \big[ e^{\beta t_{i+1}} |\Delta Y^{x,m+1}_{t_{i+1}-}|^2 \big]
		+ \frac{L_2}{\lambda_2} \| \Delta Y^{x,m}\|_{\beta, t}^2 
		+ \frac{L_2}{\lambda_2} \| \Delta Z^{x,m}\|_{\beta, t}^2. \nonumber
	\end{align}
 
	$\mathrm{iii)}$
	Let us now choose
	$1>\rho = \rho_0 > \rho_1 > \cdots > \rho_{N_h} > 0$ such that
	\be\label{eq: 1+m le }
	(m+1)e^{\beta h} \le  \big( \frac{\rho_{k}}{\rho_{k+1}} \big)^{m+1}\;\mbox{for all $m \ge 0$.}
	\ee
	\vspace{1mm}

	For $i = N_h -1$, we have $\Delta Y^{x,m}_{t_{i+1}-} = 0$ for all $m \ge 1$.
	Choosing $\lambda_1, \lambda_2$ and $\beta > 0$ in \eqref{eq:contract_beta} such that
	$$
		\frac{L_2}{\lambda_2} \frac{1}{\beta - (2L_1 + \lambda_1 L_1 + \lambda_2 L_2)} \le \rho_{N_{h}}
		\;\;\mbox{and}\;\;
		\frac{L_2}{\lambda_2} \frac{1}{2 - L_1/\lambda_1} \le \rho_{N_{h}},
	$$
	it follows from \eqref{eq:contract_beta} that, for $t\in [t_{N_h-1},t_{N_h})$,
	$
		\| \Delta Y^{x,m+1} \|_{\beta, t}^2 + \| \Delta Z^{x,m+1} \|_{\beta, t}^2
		\le
		C( \rho_i)^{m+1},
	$ for $m\ge 0$, 
	where 
	$$
	C := \esssup \sup_{(s,x')\in [0,T]\x \R^{d}} e^{\beta T}\big( | \Delta Y^{x',0}_{s} |^2 + \| \Delta Z^{x',0}_{s} \|^{2}\big)<\infty,
	$$
	and  then, by  \eqref{eq:contract_beta}  again, 
	\b*
		  |\Delta Y^{x,m+1}_{t}|^2 \le  C (\rho_i)^{m+1},
		~~\mbox{for}~t\in [t_{i},t_{i+1}), \;i = N_h -1,\;m\ge 0.
	\e*
	Recalling Remark \ref{rem: facelift non exp}, this shows that 
		\be \label{eq:Conv_DeltaY_i}
		  |\bar Y^{x,m+1}_{t}-\bar Y^{x}_{t}|^2 \le  C_{i} (\rho_i)^{m+1},
		~~\mbox{for}~t\in [t_{i},t_{i+1}),\;i = N_{h} -1,\;m\ge 0,
	\ee
	in which 
	$$
	C_{N_{h}-1}:=C.
	$$
	 
	Assume now that \eqref{eq:Conv_DeltaY_i} holds true for $i+1\le N_{h}$ and some given $C_{i+1}>0$.
	Recall that $\rho_{i}\ge \rho_{N_{h}}$. Applying \eqref{eq:contract_beta} with the above choice of  $\lambda_1, \lambda_2$ and $\beta$, we obtain
	\begin{align*}
		\| \Delta Y^{x,m+1} \|_{\beta, t}^2 + \| \Delta Z^{x,m+1} \|_{\beta, t}^2
		&\le
		e^{\beta h} C_{i+1}  (\rho_{i+1})^{m+1} \\
		&+ \rho_{i} \big(\| \Delta Y^{x,m} \|_{\beta, t}^2 + \| \Delta Z^{x,m} \|_{\beta, t}^2 \big),
	\end{align*}
	which, by \eqref{eq: 1+m le } and the fact that $\rho_{i}<1$, induces that
	\begin{align}
		\| \Delta Y^{x,m+1} \|_{\beta, t}^2 + \| \Delta Z^{x,m+1} \|_{\beta, t}^2
		&\le
		(m+1) e^{\beta h}C_{i+1} (\rho_{i+1})^{m+1}\nonumber \\
		&+ ( \rho_{i} )^{m+1}  \big(\| \Delta Y^{x,0} \|_{\beta, t}^2 + \| \Delta Z^{x,0} \|_{\beta, t}^2 \big)
		\nonumber\\
		&\le  C'_{i}\;(\rho_i)^{m+1} \label{eq: estime L2 dif YZ}
	\end{align}
	where 
	$$
	C'_{i}:=C_{i+1} + C. 
	$$
	Let us further choose $\lambda_{2}>0$ such that $L_{2}/\lambda_{2}\le \rho_{N_{h}}$, and recall that $\rho_{i}\ge  \rho_{N_{h}}$. 
	Then, using again \eqref{eq:contract_beta}, \eqref{eq: 1+m le }, \eqref{eq:Conv_DeltaY_i} applied to $i+1$, we obtain, for $t\in [t_{i},t_{i+1})$,
	\begin{align*}
		  |\Delta Y^{x,m+1}_{t}|^2  
		&\le e^{\beta h} C_{i+1}( \rho_{i+1})^{m+1} + C_{i}'  (\rho_i)^{m+1}\le 
		C_{i}(\rho_i)^{m+1}, 
	\end{align*}
	so that it follows from Remark \ref{rem: facelift non exp}  that
	\begin{align}
		  |\bar Y^{x,m+1}_{t}-\bar Y^{x}_{t}|^2& \le 
		C_{i}(\rho_i)^{m+1},\label{eq: estime sup dif Y}
	\end{align}
	
	where
	$$
	C_{i}:=  e^{\beta h}C_{i+1} + C_{i}'.
	$$
	Since $(\Yb,\Zb)=(\Yb^{X_{t_{i}}}, \Zb^{X_{t_{i}}})$ and $(\Yb^{m},\Zb^{m})=(\Yb^{X_{t_{i}},m},\Zb^{X_{t_{i}},m})$ on each $[t_{i},t_{i+1})$, this concludes the proof.
	\qed

\section{A branching diffusion representation for $\Yb^m$}
\label{sec: repre}

	We now explain how the  solution of   \eqref{eq:def_Ym} on $[t_i, t_{i+1})$ can be represented by means of a branching diffusion system. We slightly adapt the arguments of \cite{HOTTW}. 
	
	Let  us  consider an element  {$(p_\ell)_{\ell \in L} \subset \R_+$ such that $\sum_{\ell \in L} p_{\ell} = 1$},   set $K_{n}:=\{(1,k_{2},\ldots,k_{n}): (k_{2},\ldots,k_{n})\in \{1,\ldots,(q_{\circ}+1)L_{\circ}\}^{n-1}\}$ for $n\ge 1$, and $K:=\cup_{n\ge 1} K_{n}$. Let $(W^{k})_{k\in K }$ be a sequence of independent $d$-dimensional Brownian motions, $(\xi_{k}=(\xi_{k,q})_{0\le q\le q_{\circ}})_{k\in  K }$ and $(\delta_{k})_{k\in  K}$ be two sequences of independent random variables, such that 
$$
\P[\xi_{k}=\ell]=p_{\ell}, \,\;\ell \in L, k\in  K , 
$$
and 
\b*
\bar F(t):=\P[\delta_{k}>t]=\int_{t}^{\infty} \rho(s)ds,\; t\ge 0,\;k\in K , 
\e*
for some continuous strictly positive map $\rho: \R_+ \to \R_+$. 
We assume that 
\b*
\mbox{$(W^{k})_{k\in K }$, $(\xi_{k})_{k\in K }$, $(\delta_{k})_{k\in K }$ and $W$ are independent.}
\e*
Given the above, we construct particles $X^{(k)}$ that have the dynamics \reff{eq: Diffusion} up to a killing time $T_{k}$ at which they split in $\|\xi_{k}\|_{1}:=\xi_{k,0}+\cdots \xi_{k,q_{\circ}}$ different (conditionally) independent particles with dynamics   \reff{eq: Diffusion} up to their own killing time. The construction is done as follows. 
First, we set $T_{(1)}:=\delta_{1}$, and, given $k=(1,k_{2},\ldots,k_{n}) \in K_{n}$ with $n\ge 2$, we let $T_{k}:=\delta_{k}+T_{k-}$ in which $k-:=(1,k_{2},\ldots,k_{n-1}) \in K_{n-1}$. We can then define the Brownian particles $(W^{(k)})_{k\in K}$ by using the following induction: we first set 
$$
W^{((1))}:=W^{1}\1_{[0,T_{(1)}]}\;,\;\Kc^{1}_{t}:=\{(1)\}\1_{[0,T_{(1)}]}(t)+\emptyset  \1_{[0,T_{(1)}]^{c}}(t), \;t\ge 0,
$$
then, given $n\ge 2$ and  $k\in \bar \Kc^{n-1}_{T}:=\cup_{t\le T}\Kc^{n-1}_{t}$, we let 
\b* 
W^{(k\oplus j)}:=\left(W^{(k)}_{\cdot\wedge T_{k}} +W^{{k\oplus j}}_{\cdot \vee  T_{k}}-W^{{k\oplus j}}_{T_{k}}\right)\1_{[0,T_{k\oplus j}]}, \; 1\le j\le \|\xi_{k}\|_{1},
\e*
in which we use the notation $$(1,k_{1},\ldots,k_{n-1})\oplus j:= (1,k_{1},\ldots,k_{n-1},j),$$ 
and
$$
	\bar \Kc^n_{t}:=\{k\oplus j: k\in \bar \Kc^{n-1}_{T}, 1\le j\le \|\xi_{k}\|_{1} ~\mbox{s.t.}~ t \in (0, T_{k \oplus j}] \},
	\;\;
	\bar \Kc_{t} := \cup_{n \ge 1} \bar \Kc^n_{t},
$$
\b*
	\Kc^n_t := \{ k\oplus j: k\in \bar \Kc^{n-1}_{T}, 1\le j\le \|\xi_{k}\|_{1} ~\mbox{s.t.}~ t\in (T_{k},T_{k\oplus j}]\},
	\;\;
	\Kc_t := \cup_{n\ge 1} \Kc^n_t.
\e*

	Now observe that the solution $X^{x}$ of \reff{eq: Diffusion} on $[0,T]$ with initial condition $X^{x}_{0}=x\in \R^{d}$ can be identified in law on the canonical space as a process of the form $\Phi[x](\cdot,W)$ in which the deterministic map $(x,s,\omega)\mapsto \Phi[x](s,\omega)$ is $\Bc(\R^{d}) \otimes \Pc$-measurable, where $\Pc$ is the predictable $\sigma$-filed on $ [0,T]\x \Omega$. We then define the corresponding particles $(X^{x,(k)})_{k\in K}$ by $X^{x,(k)}:=\Phi[x](\cdot,W^{(k)})$.
	Moreover, we define the $d\x d$-dimensional matrix valued tangent process $\nabla X^{x,(k)}$ defined on $[\![T_{k-},T_{k}]\!]$ by
	\begin{align} \label{eq:TanProc}
		\nabla X^{x,(k)}_{T_{k-}} &= I_d\\
		 d\nabla X^{x,(k)}_t &= D \mu(X^{x,(k)}_t) \nabla X^{x,(k)}_t dt + \sum_{i=1}^d D \sigma_i(X^{x,(k)}_t) \nabla X^{x,(k)}_t dW^{(k),i}_t,\nonumber
	\end{align}
	where $I_d$ denotes the $d \x d$-dimensional identity matrix, 
	and $\sigma_i$ denotes the $i$-th column of matrix $\sigma$.
\vs2

Finally, we give a mark $0$ to the initial particle $(1)$,
	and, for every particle $k$, knowing $\xi_k = (\xi_{k,0}, \cdots, \xi_{k, q_{\circ}})$,
	we consider its offspring particles $(k \oplus j)_{ j = 1, \cdots, \|\xi_k\|_{1}}$ and give the first $\xi_{k,0}$ particles the mark $0$, the next $\xi_{k,1}$ particles the mark $1$, the next $\xi_{k,2}$ particles the mark $2$, etc.
	Thus, every particles  $k$ carries a mark $\theta_{k}$ taking values in $0$ to $q_{\circ}$. 
\vspace{2mm}

Given the above construction, we can now
 provide the branching process based representation of $(\bar Y^{m})_{m\ge 0}$. We assume here that $(\bar u^{m-1},\bar v^{m-1})$ defined in \eqref{eq: def bar um bar vm} are given for some $m\ge 1$, recall that 
$(\bar u^{0},\bar v^{0})=({\rm y},D{\rm y})$ by construction. 
We  set $\tilde u^{m}(T,\cdot)=g$. Then, for $i=N_{h}-1,\cdots,0$, we define $(\tilde u^{m},\tilde v^{m})$ on each interval $[t_{i},t_{i+1})$ recursively by 
\begin{align}
	\tilde u^{m}(t,x)&:=\E\left[U^{m}_{t,x}\right]\1_{t\ne t_{i}}+ \Rc[\E\left[U^{m}_{t,\cdot}\right]](x)\1_{t= t_{i}}\nonumber
	\\
	\tilde v^m(t,x)& := \E \big[ V^m_{t,x} \big]\sigma(x), \label{eq:repres_vm}
\end{align}
for $(t,x)\in [t_{i},t_{i+1})\x \R^{d}$, in which 
\b*
	U^{m}_{t,x}
	&:=& 
	\Big[\prod_{k \in   \Kc_{\tip-t}}G^{m}_{t,x}(k) \Wc_{t,x}(k) \Big] \Big[\prod_{k \in\bar{\Kc}_{\tip-t}\setminus\Kc_{\tip-t}}A^{m}_{t,x}(k) \Wc_{t,x}(k)\Big] \\
	V^m_{t,x}
	&:=&
	U^m_{t,x} ~\frac{1}{T_{(1)}} \int_0^{T_{(1)}}  \big[ \sigma^{-1}(X^{x,(1)}_s) \nabla X^{x,(1)}_s \big]^{\top} dW^{(1)}_r
\e*
	where
$$
	G^{m}_{t,x}(k)
	~:=~
	\frac{\tilde u^{m} \big(t_{i+1}, X_{\tip-t}^{x,(k)} \big) - \tilde u^{m} \big(t_{i+1}, X_{T_{k-}}^{x,(k-)} \big) \1_{\{\theta(k) \neq 0\} }}
	{\bar{F}({\tip-t}-T_{k-})},
$$
$$
	A^{m}_{t,x}(k)
	:= 
	\frac{\sum_{j=1}^{j_{\circ}}c_{j,\xi_k} \big(X_{T_k}^{x,(k)} \big)
			\vp_{j} \big( (\bar u^{m-1}, \bar  v^{m-1}) (t+T_{k},X_{T_k}^{x,(k)}) \big)
	}{p_{\xi_k}\,\rho(\delta_k)},
$$
$$
	\Wc_{t,x}(k) 
	:=
	\1_{\{\theta_k = 0\}}
	+ \frac{ \1_{\{\theta_k \neq 0\}}}{T_k - T_{k-}} 
	b_{\theta_k}(X^{x,(k)}_{T_k}) \cdot \int_{T_{k-}}^{T_k} \big[ \sigma^{-1}(X^{x,(k)}_s) \nabla X^{x,(k)}_s \big]^{\top} dW^{(k)}_r.
$$
 Compare with \cite[(3.4) and (3.10)]{HOTTW}.

	\vspace{0.5em}

	The next proposition shows that $(\tilde u^{m},\tilde v^{m})$ actually coincides with $(\bar u^{m},\bar v^{m})$
	in \eqref{eq: def bar um bar vm},
	a result that follows essentially from \cite{HOTTW}.
	Nevertheless, to be more precise on the square integrability of $U^m_{t,x}$ and $V^m_{t,x}$, 
	one will fix a special density function $\rho$ as well as  probability weights $(p_{\ell})_{\ell \in L}$.
	Recall again that $(\bar u^m, \bar v^m)$ are defined in Theorem \ref{thm: main} and satisfy \eqref{eq: def bar um bar vm},
	and that $(h_{\circ},M_{\circ})$ are chosen such that \eqref{eq: def Mho et ho}-\eqref{eq: def Mho et ho Z} in the Appendix hold.

\begin{Proposition}\label{prop: representation Ynm}
	Let us choose $\rho(t) = \frac{1}{3} t^{-2/3} \1_{\{t \in [0,1]\}}$ and
	$p_{\ell} = \frac{\|c_{\ell}\|_{\infty}}{\|c\|_{1,\infty}}$ with $\|c\|_{1,\infty} := \sum_{\ell \in L} \|c_{\ell}\|_{\infty}$. Let  $h'_{\circ}$ and $M_{h'_{\circ}}$ be given by \eqref{eq:def_delta} and \eqref{eq:def_M_delta}. Assume that $h\in (0,h_{\circ}\wedge h'_{\circ})$.  
	Then, 
	$$
	\E[|U^{m}_{t,x}|^2] \vee \E[\|V^{m}_{t,x}\|^2] \le   (M_{h'_{\circ}})^2,\;\mbox{ for all   $m \ge 1$ and  $(t,x) \in [0,T] \x\R^d$.}
	$$
	Moreover,  $\tilde u^{m} = \bar u^{m}$ on $[0,T]\x \R^{d}$ and $\tilde  v^{m}= \bar v^m$ on $\cup_{i\le N_{h}-1}(t_{i},t_{i+1})\x \R^{d}$.
\end{Proposition}

The proof of the above mimics the arguments of   \cite[Theorem 3.12]{HOTTW} and is postponed to the Appendix \ref{Proof of the representation formula}. 
 
\begin{Remark}
	The integrability and representation results in Proposition \ref{prop: representation Ynm}  hold true for a large class of parameters $\rho$, $(p_{\ell})_{\ell \in L}$ and $(c_{\ell})_{\ell \in L}$ (see e.g.~\cite[Section 3.2]{HOTTW} for more details).
	We restrict to a special choice of parameters in Proposition \ref{prop: representation Ynm} in order to compute explicitly the lower bound $h'_{\circ}$ for the explosion time  as well as the upper bound $M_{h'_{\circ}}$ for the variance of the estimators.
\end{Remark}

 \begin{Remark}\label{rem: numeric} a. The above scheme requires the computation of conditional expectations that involve the all path $(\tilde u^{m},\tilde v^{m})(\cdot,X^{t_{i},x})$ on each $(t_{i},t_{i+1})$, for all $x\in \R^{d}$, and therefore the use of an additional time space grid on which $\tilde u^{m}$ and $\tilde v^{m}$ are estimated by Monte-Carlo. For $t_{i}<t<t_{i+1}$, the precision does not need to be important because the corresponding values are only used for the localization procedure, see the discussion just after Remark \ref{rem: facelift non exp}, and the grid does not need to be fine. The space grid should be finer at the times in $\T$ because each $\tilde u^{m}(t_{i},\cdot)$ is used per se as a terminal condition on $(t_{i},t_{i+1}]$, and not only for the localization of the polynomial.   This corresponds to the Method B in \cite[Section 3]{BouchardTanZou}. One can also consider a very fine time grid  $\T$ and avoid the use of a sub-grid. This is Method A in \cite[Section 3]{BouchardTanZou}. The numerical tests performed in \cite{BouchardTanZou} suggest that Method A is more efficient. 
 
b. From a numerical viewpoint, $\tilde v^{m}$ can also be estimated by using a finite difference scheme based on the estimation of $\tilde u^{m}$. It seems to be indeed more stable. 

c. Obviously, one can also adapt the algorithm to make profit of the ghost particle or of the  normalization techniques described in \cite{WBranchVar}, which seem to  reduce efficiently the variance of the estimations, even when $\rho$ is the exponential distribution.   

 \end{Remark}

\section{Numerical  example}
\label{sec:NumericalExamples}
This section is dedicated to a simple example in dimension one showing the efficiency of the proposed methodology. We use the Method A of \cite[Section 3]{BouchardTanZou} together with the normalization technique described in \cite{WBranchVar} and a finite difference scheme to compute $\tilde v^{m}$.
	In particular, this normalization technique allows us to take $\rho$ as an exponential density rather than that in Proposition \ref{prop: representation Ynm}.
	See Remark \ref{rem: numeric}.
 \vs2
 
We consider the SDE with coefficients
\begin{align*}
\mu(x) = -0.5(x+0.2)\;\mbox{ and }\; \sigma(x) =0.1 \left(1+ 0.5 ((x\vee 0)\wedge 1) \right).
\end{align*}
The maturity is   $T=1$ and the non linearity in \eqref{eq:BSDE} is taken as
\begin{flalign}
f(t,x,y,z)= & \hat f(y,z) + \frac{1}{2} e^{\frac{t-T}{2}} \left( \cos(x) \frac{\sigma(x)^2}{2}  - \frac{1}{2} (1+\cos(x)) + 
  \mu(x) \sin(x) \right)   \nonumber \\
&-  \frac{2}{4+ |\sin(x) (1+\cos(x))| e^{t-T}} \nonumber 
\end{flalign}
where
\begin{flalign}
\hat f(y,z) = &  \frac{1}{2(1+|yz|)} .
\label{eq:fExample}
\end{flalign}
It is complemented by the choice of the terminal condition $g(x)= \frac{1+ \cos(x)}{2}$, so that an analytic solution is available :
\begin{flalign*}
u(t,x) =  \frac{1+ \cos(x)}{2} e^{\frac{t-T}{2}}.
\end{flalign*}
We use the algorithm to compute an estimation of $u(0,\cdot)$ on $\Xbf:=[-1.,1.]$.

\begin{figure}[H]
\begin{minipage}[b]{0.5\linewidth}
  \centering
 \includegraphics[width=0.9\textwidth]{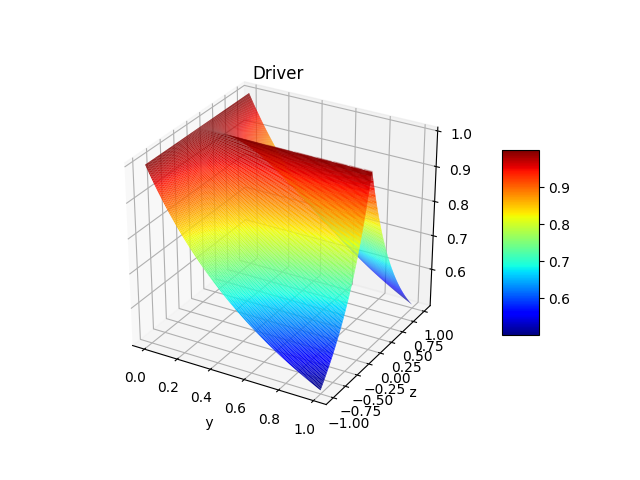}
 \caption*{\small The driver $\hat f$}
 \end{minipage}
\begin{minipage}[b]{0.5\linewidth}
  \centering
 \includegraphics[width=0.9\textwidth]{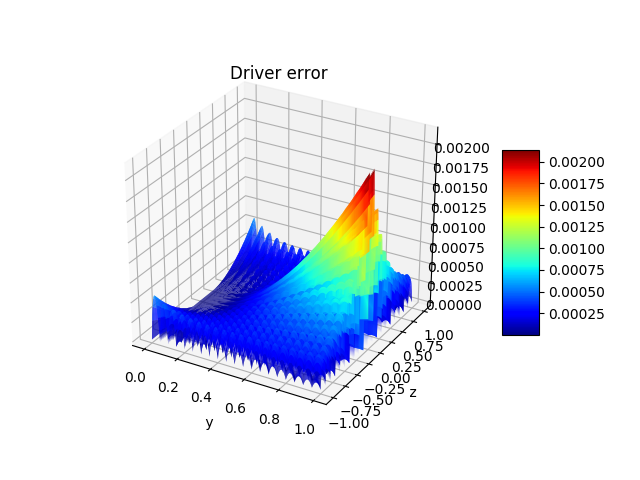}
 \caption*{\small Error on the driver due to the linear spline representation.}
 \end{minipage}
\caption{The driver $\hat f$ and its linear  spline representation error for $20 \times 20$ splines. } 
\label{fig:driver}
\end{figure}

To construct our local polynomials approximation of $\hat f$, we use a linear spline interpolation in each direction, obtained by tensorization, and
leading to a local approximation on the basis $1$, $y$, $z$, $yz$ on each mesh of a domain $[0,1] \times [-1,1]$.
Figure \ref{fig:driver} displays  the function $\hat f$ and the error obtained by a discretization of $20 \times 20$ meshes.

The parameters affecting the convergence of the algorithm are:
\begin{itemize}
\item The couple of meshes $(n_y, n_z)$ used in the spline representation  of \eqref{eq:fExample}, where $n_y$ (resp.~$n_{z}$) is the number of linear spline meshes for the $y$ (resp.~$z$) discretization.
\item The number of time steps  $N_h$.
\item The grid and the interpolation used on $\Xbf$ at   $t_0=0$ and for all dates $t_i$, $i=1,\ldots,N_h$. Note that  the size of the grid has to be adapted to the value of $T$, because of the diffusive feature of \eqref{eq: Diffusion}. All interpolations are achieved with the StOpt library (see \cite{GLWa,GLWa2}) using a modified quadratic interpolator as in \cite{warin}. In the following, $\Delta x$ denotes the mesh of the space discretization.
\item The time step  is set to $0.002$ and we use an Euler scheme to approximate \eqref{eq: Diffusion}.
\item The accuracy of the  estimation of the expectations appearing in our algorithm. We compute the empirical standard deviation $\theta$ associated to each Monte Carlo estimation of 
the expectation in \reff{eq:repres_vm}. We try to  fix the number $\hat M$ of samples such that  ${\theta}/{\sqrt{\hat M}}$ does not exceed a certain level, fixed at $0.000125$, at each point  of our grid. We cap  this number of simulations at  $510^{5}$.
\item  The intensity, set to $0.4$, of the exponential distribution  used to define the random variables $(\delta_k)_{k\in K}$.
\end{itemize} 
Finally, we take $M=1$ in the definition of $\Rc$.
\vs2

We only perform one Picard iteration with initial prior $(\tilde u^{0},\tilde v^{0})=(g,Dg\sigma)$. 
\vs2

\begin{figure}[H]
 \centering
\includegraphics[width=0.45\textwidth]{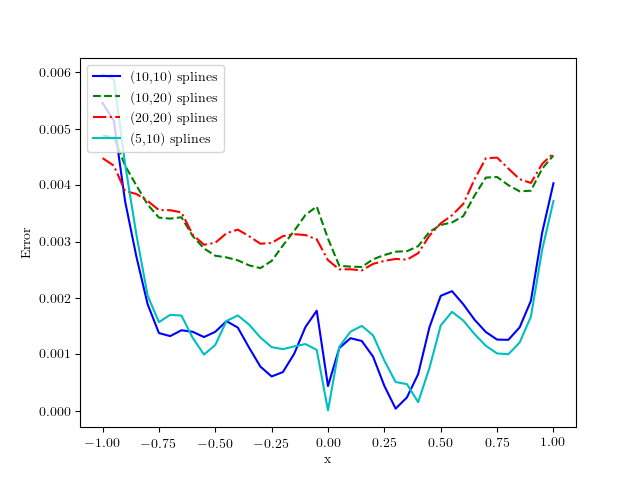}
\includegraphics[width=0.45\textwidth]{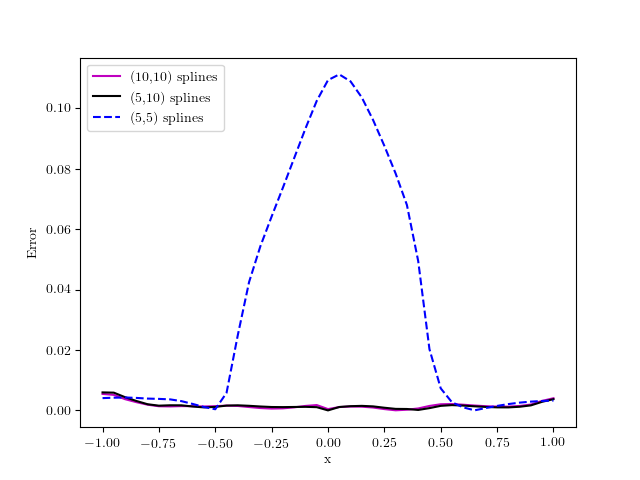}
 \caption{Error plot depending on the couple $(n_y,n_z)$ for $20$ time steps,  a space discretization  $ \Delta x= 0.1$.}
\label{fig:errorOnSpline}
\end{figure}

\begin{figure}[H]
\centering
\includegraphics[width=0.45\textwidth]{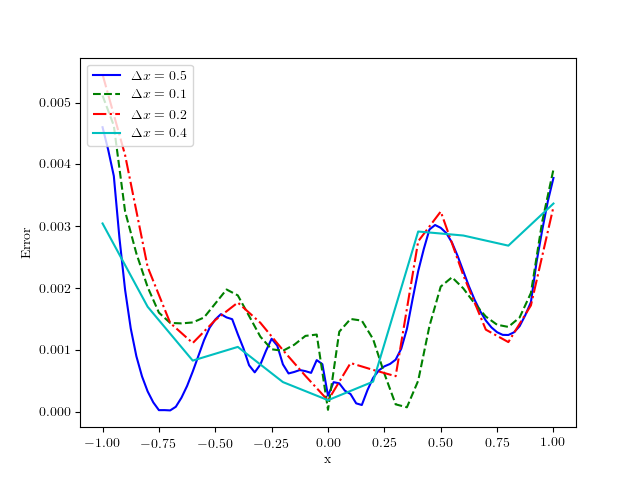}
 \caption{Error plot depending on $\Delta x$ for $(n_y,n_z) = (20,10)$, $N_h=20$.}
\label{fig:errorOnInterpol}
\end{figure}

\begin{figure}[H]
 \centering
\includegraphics[width=0.45\textwidth]{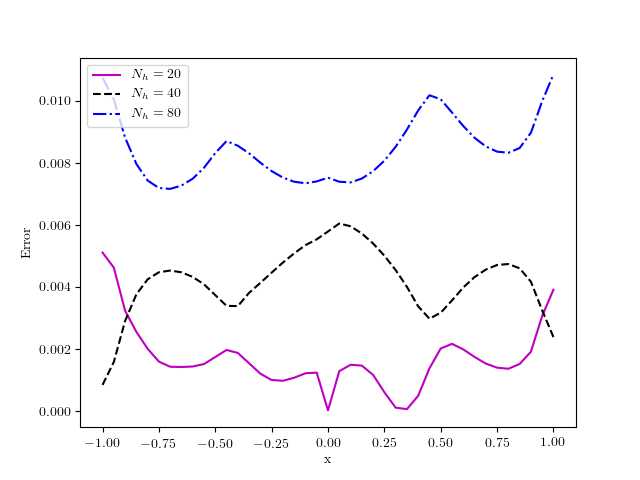}
\includegraphics[width=0.45\textwidth]{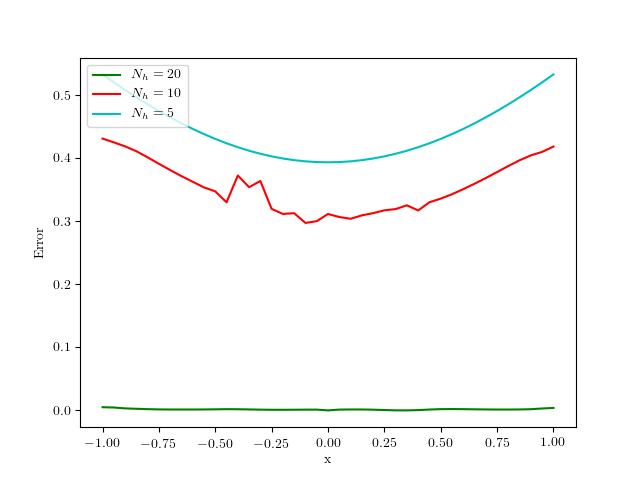}
 \caption{Error plot depending on $N_h$  for  $(n_y,n_z) =(20,10)$, a space discretization  $ \Delta x= 0.1$}
\label{fig:errorOnTStep}
\end{figure}

On the different figures below, we plot the errors obtained on $\Xbf$ for different values of  $N_h$,  $(n_y, n_z)$ and  $\Delta x$.
We first    use $20$ time steps and an interpolation step of $0.1$ In  figure  \ref{fig:errorOnSpline},  we display  the error as a function of the number of spline meshes.
We provide two plots:
\begin{itemize}
\item On the left-hand side,  $n_y$ varies above 5 and $n_z$ varies above 10,
\item On the right-hand side, we add  $ (n_y,n_z)=(5,5)$. It leads to   a maximal error of $0.11$, showing that a quite good accuracy in the spline representation in $z$ is necessary.
\end{itemize}
In figure \ref{fig:errorOnInterpol}, we plot the error obtained with $(n_y,n_z) = (20,10)$ and a number of time steps equal to $N_h=20$, for different values of $\Delta x$: the results are remarkably stable
with the interpolation space discretization.
In Figure  \ref{fig:errorOnTStep}, we finally let the number of time steps  $N_h$ vary. Once again we give two plots:
\begin{itemize}
\item one with  $N_h$ above or equal to $20$,
\item one with  small values of  $N_h$.
\end{itemize}
The results clearly show that the algorithm produces bad results when $N_h$ is too small: the time steps are too large for the branching method. In this case, it  exhibits a large variance. When $N_h$ is too large, then interpolation errors  propagate leading also to a deterioration of our estimations.
Numerically, it can be checked that the face-lifting procedure is in most of the cases useless when only one Picard iteration is used:
\begin{itemize}
\item When the variance of the branching scheme is small, the face-lifting and truncation procedure has no effect,
\item When the variance becomes too large, the face-lifting procedure is  regularizing the solution and this permits to reduce the error due to our interpolations.
\end{itemize}

\begin{figure}[H]
 \centering
\includegraphics[width=0.45\textwidth]{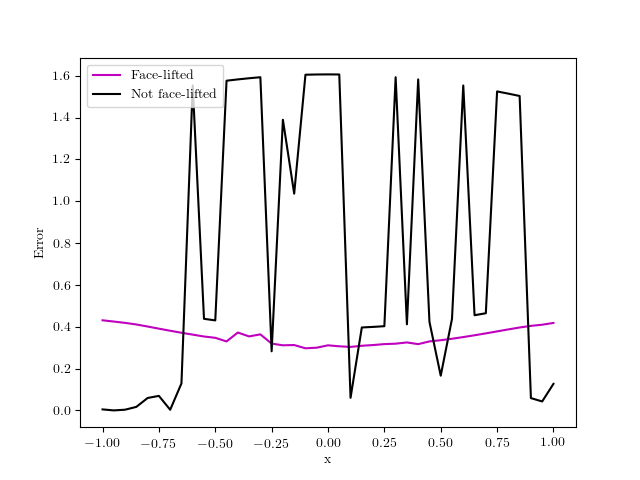}
\caption{Error plot  for  $N_h=10$  for  $(n_y,n_z) =(20,10)$, a space discretization  $ \Delta x= 0.1$ with and without face-lifting.}
\label{fig:faceLift}
\end{figure}

In figure \ref{fig:faceLift}, we provide the estimation with  and without face-lifting, obtained with $N_h=10$,  $(n_y,n_z) =(20,10)$ and a space discretization  $ \Delta x= 0.1$.

For $N_h=20$ the computational time is equal to  $320$ seconds on a  $3$ year old laptop.

\appendix
\section{Appendix}
\subsection{A priori estimates for the Picard iteration scheme}\label{sec: A priori estimates for the Picard iteration scheme}
\def\ut{\tilde u}

In this section, we let $\nabla X$ be the tangent process associated to  $X$ on $[0,T]$ by 
\begin{align*} 
		\nabla X_{0} &= I_d\;,\;
		 d\nabla X_t = D \mu(X_t) \nabla X_t dt + \sum_{i=1}^d D \sigma_i(X_t) \nabla X_t dW^{i}_t,\nonumber
	\end{align*}
and we define 
$$
N_{s}^{t}:=\left(\int_{t}^{s} \frac{1}{s-t}  (\sigma(X_{r})^{-1} \nabla X_{r})^{\top} dW_{r}\right)^{\top} \nabla X_{t}^{-1} 
$$
for $t\le s\le T$. Standard estimates lead to 
\begin{align}\label{eq: borne Nst}
\E_{t}[\|N_{s}^{t}\| ]\le  C_{\mu,\sigma}(s-t)^{-\frac12} \;\mbox{ for } t\le s \le T, 
\end{align}
for some $C_{\mu,\sigma}>0$ that only depends on  $\|\mu\|_{\infty}$, $\|\sigma\|_{\infty}$, $\|D\mu\|_{\infty}$, $\|D\sigma\|_{\infty}$, $a_{0}$ in \eqref{eq: unif elliptic} and $T$. In particular, it does not depend on $\|X_{0}\|$. 
Up to changing this constant, we may assume that 
\begin{align}\label{eq: borne nabla X}
\E_{t}[\|\nabla X_{s}\nabla X_{t}^{-1}\|]\le e^{C_{\mu,\sigma} (s-t)} \;\mbox{ for } t \le s\le T. 
\end{align}

Set
$$
\Dc_{M_{h_{\circ}}}:=\R^{d}\x \{(y,z)\in \R\x \R^{d}: |y|\vee \|z^{\top} \sigma^{-1}\|\le M_{h_{\circ}}\}^{2}.
$$
and let $M_{h_{\circ}}\ge M$ and $h_{\circ}\in (0,T]$ be such that 
\begin{align}
 M_{h_{\circ}}&\ge  M+ \sup_{\Dc_{M_{h_{\circ}}}}|f_{\circ}|\;h_{\circ}.\label{eq: def Mho et ho}
\end{align}
and
\begin{align}
 M_{h_{\circ}}&\ge    M  e^{C_{\mu,\sigma} h_{\circ}}+C_{\mu,\sigma}\; \sup_{\Dc_{M_{h_{\circ}}}}|f_{\circ}|\;(h_{\circ})^{\frac12} \label{eq: def Mho et ho Z}
\end{align}
The existence of $h_{\circ}$ and $M_{h_{\circ}}$ follows from  \eqref{eq: f loc poly} and \eqref{eq: borne a vp}. Note that they do not depend on $X_{0}$.

\begin{Proposition}\label{prop: bound tilde Y tilde Z} Let $\tilde g:\R^{d}\mapsto \R$ be bounded by $M$ and $M$-Lipschitz. Fix $T'\in [0,h_{\circ}]$.
Let $(\tilde U,\tilde V ): [0,T']\x \R^{d}\mapsto \R\x\R^{d}$ be measurable such that $|\tilde U|\vee \|\tilde V^{\top}\sigma^{-1}\|\le M_{h_{\circ}}$.  Then, there exists a unique bounded solution on  $[0,T']$ to 
\begin{align}\label{eq: bsde tilde Y}
\Yt_{\cdot}
	=
	\tilde g({X}_{T'})
	+\!
	\int_{\cdot}^{T'} f_{{\circ}}( {X}_s,\Yt_s, \Zt_s, (\tilde U,\tilde V)(s,X_{s}))\,ds
	\!-\!\int_{\cdot}^{T'} \Zt_s^{\top}\,dW_s.
\end{align}
It satisfies 
\begin{align}
|\Yt|\vee\|\Zt^{\top}\sigma(X)^{-1}\| \le M_{h_{\circ}} \;\mbox{on } [0,T'].
\end{align}
Moreover, there exists a bounded continuous map $(\tilde u,\tilde v) : [0,T']\x \R^{d}\mapsto \R\x\R^{d}$ such that $\Yt=\tilde u(\cdot, X)$ on $[0,T']$ $\P$-a.s.~and $\Zt=\tilde v(\cdot, X)$ $dt\x d\P$-a.e.~on $[0,T']\x \Omega$. It satisfies $\tilde v^{\top}=D\tilde u \sigma$ on  $[0,T')$.
\end{Proposition}

\proof   We construct the required solution by using Picard iterations. We set $(\Yt^{0},\Zt^{0})=({\rm y}, D{\rm y})(\cdot,X)$, and define recursively on $[0,T']$ the couple $(\Yt^{n+1},\Zt^{n+1})$ as the unique solution of 
\begin{align*}
\Yt^{n+1}_{\cdot}
	=
	\tilde g({X}_{T'})
	+\!
	\int_{\cdot}^{T'} f_{{\circ}}( {X}_s,\Yt^{n}_s, \Zt^{n}_s, (\tilde U,\tilde V)(s,X_{s}))\,ds
	\!-\!\int_{\cdot}^{T'} (\Zt^{n+1}_s)^{\top}\,dW_s, 
\end{align*}
whenever it is well-defined. It is the case for $n=0$. We now assume that $(\Yt^{n},\Zt^{n})$ is well-defined and such that $|\Yt^{n}|\vee\|\sigma(X)^{-1}\Zt^{n}\| \le M_{h_{\circ}}$ on $[0,T']$ for some $n\ge 0$. 
Then,  
\begin{align*}
|\Yt^{n+1}_{\cdot}|\le \|\tilde g\|_{\infty}+\sup_{\Dc_{M_{h_{\circ}}}}|f_{\circ}| \;h_{\circ}  \le M_{h_{\circ}},
\end{align*}
in which we used \eqref{eq: def Mho et ho} for the second inequality. On the other hand, up to using a mollifying argument, one can assume that $\tilde g$ is $C^{1}_{b}$ and that $(U,V)$ is Lipschitz. Then, it follows from the same arguments as in \cite[Theorem 3.1, Theorem 4.2]{MaZhangRefl} that $\Zt^{n+1}$ admits the representation 
\begin{align*}
( \Zt^{n+1}_{t})^{\top}=&\E_{t}\left[D\tilde g(X_{T'})\nabla X_{T'}\nabla X_{t}^{-1}\right]\sigma(X_{t})\\
&+\E_{t}\left[\int_{t}^{T'} f_{{\circ}}( {X}_s,\Yt^{n}_s, \Zt^{n}_s, (U,V)(s,X_{s})) N_{s}^{t} ds \right]\sigma(X_{t}).
\end{align*}
By combining the above together with \eqref{eq: borne Nst} and \eqref{eq: borne nabla X}, we obtain that 
\begin{align*}
 \|(\Zt^{n+1}_{t})^{\top}\sigma(X_{t})^{-1}\|  &\le    M  e^{C_{\mu,\sigma} (T'-t)}+C_{\mu,\sigma}\; \sup_{\Dc_{M_{h_{\circ}}}}|f_{\circ}|\;(h_{\circ})^{\frac12}   \le M_{h_{\circ}},
\end{align*}
in which we used \eqref{eq: def Mho et ho Z} for the second inequality.  
	The above proves that the sequence $(\Yt^{n},\Zt^{n})_{n\ge 0}$ is uniformly bounded on $[0,T']$. 
	Therefore, we can consider $f_{\circ}$ as a Lipschitz generator, and hence $(\Yt^{n},\Zt^{n})_{n\ge 0}$
	is in fact a Picard iteration that converges to a solution of \eqref{eq: bsde tilde Y} with the same bound.

The existence of the   maps 
$\tilde u$ and $\tilde v $ such that $\tilde v^{\top}=D\tilde u\sigma$ follows from \cite[Theorem 3.1]{MaZhangRefl} applied to \eqref{eq: bsde tilde Y} when $(\tilde g,  f_{\circ}(\cdot,y,z,(U,V)(t,\cdot)))$ is $C^{1}_{b}$, uniformly in  $(t,y,z)\in [0,T]\x\R^{d+1}$. The representation result of \cite[Theorem 4.2]{MaZhangRefl} combined with a simple approximation argument, see e.g.~\cite[(ii) of the proof of Proposition 3.2]{FLLT},  then shows that the same holds on  $[0,T')$ under our conditions. 
\ep

\subsection{Proof of the representation formula}\label{Proof of the representation formula}

	We adapt the proof of    \cite[Theorem 3.12]{HOTTW} to our context. We proceed by induction.  In all this section, we fix 
	$$
	(t,x)\in [t_{i},t_{i+1})\x \R^{d}, 
	$$
	and assume that the result of Proposition \ref{prop: representation Ynm} holds up to rank $m-1\ge 0$ on $[0,T]$ (with the convention $U^{0}_{\cdot}=y$, $V^{0}_{\cdot}:=Dy$),  and  up to rank $m$ on   $[t_{i+1},T]$. In particular,  
	 we assume that $\tilde u^{m}(t_{i+1},\cdot) = \bar u^{m}(t_{i+1},\cdot)$.
	 
	We fix $\varrho=4$ and define  
	\begin{align*}
		C_{1,\varrho} 
		&:= 
		 M^{\varrho} \vee \sup_{t \le s,~ x\in \R^d, ~q=1, \cdots, q_{\circ}, ~\|\xi\|\le M}
		\E \Big[ \Big| \big (\xi \cdot (X_s^{t,x} -x) \big) \big( b_q(x) \cdot \overline \Wc_{t,x,s} \big) \Big|^\varrho \Big],
	\\
		C_{2,\varrho} 
		&:=
		\sup_{t \le s\le t_{i+1},~ x\in\R^d, ~q=1, \cdots,q_{\circ}} 
		\E \Big[ \big| \sqrt{s-t}~ b_q(x) \cdot  \overline \Wc_{t,x,s} \big|^\varrho \Big]
			\end{align*}
	where 
		$$
		\overline \Wc_{t,x,s}
		:=
		\frac{1}{s-t} \int_t^s \big[ \sigma^{-1}(X^{t,x}_r) \nabla X^{t,x}_r \big]^{\top} dW_r
	$$
	in which   $\nabla X^{t,x}$ is the tangent process of $X^{t,x}$ with initial condition $I_{d}$ at $t$. We then set 
	\begin{align*}
		\hat C_{1,\varrho} := \frac{C_{1,\varrho}}{\bar F(T)^{\varrho-1}},
		\;\;\;\;\;
		\hat C_{2,\varrho} := C_{2,\varrho} \;j_{\circ} \sup_{j\le j_{\circ},\ell \in L, t\in (0, h_{\circ}]} \Big( \frac{\|c_{j,\ell}\|_{\infty}}{p_{\ell}} \frac{t^{-\varrho/(2(\varrho-1))}}{\rho(t)} \Big)^{\varrho-1}.
	\end{align*}
 Since $\bar F$ is non-increasing and $\tilde u^{m}(t_{i+1},\cdot)=\bar u^{m}(t_{i+1},\cdot)$ is bounded by $M$ and $M$-Lipschitz, direct computations imply that 
 \begin{align}
 \E[|U^{m}_{t,x}|^{\varrho}]\vee\E\|V^{m}_{t,x}\|^{\varrho}] \le \E\Big[
		\Big( \prod_{k \in \Kc_t} \frac{\hat C_{1,\varrho}}{\bar F(t - T_{k-})} \Big)
		\Big( \prod_{k \in \bar \Kc_t \setminus \Kc_t} \frac{ \hat C_{2,\varrho}}{p_{\xi_k} \rho(\delta_k)} \Big)
		\Big].\label{eq: majo puissance 4}
 \end{align}
 
 We first estimate the right-hand side, see \eqref{eq:branch_bounded_by_Mh} below.  

	Let us denote by $C_{\rm bdg}$ the constant in the Burkholder-Davis-Gundy inequality such that $\E\big[ \sup_{0 \le t \le T} |M_t|^\varrho] \le C_{\rm bdg} \E[ (\langle M \rangle_T)^{\frac{\varrho}{2}} ]$ for any continuous martingale $M$ with $M_0 = 0$.
	Denote further 
	$$C_0 :=  (3 \x 3)^{\varrho-1} \big( 1+ (\varrho/(\varrho-1))^\varrho \big) \big( 1 + (1 + |\bar \lambda_{D\mu} T|^\varrho e^{C_Q T} ) \big),$$
	where $\bar \lambda_{D\mu}$ the largest eigenvalue of the matrix $D \mu$, $\bar \lambda_{D\sigma}$ the largest eigenvalue of the matrix $D \sigma_i$, $i\le d$,
	and $C_Q := \varrho  \bar \lambda_{D\mu} + d \varrho(\varrho-1) \bar \lambda_{D \sigma}/2$.
	Define also $\bar \lambda_{(\sigma \sigma^{\top})^{-1}}$ as the largest eigenvalue of matrix $(\sigma \sigma^{\top})^{-1}$.
	
	\begin{Lemma} \label{lemm:C12q}
		Under the Assumptions of Proposition \ref{prop: representation Ynm},  
	$$
		\hat C_{1,\varrho}\le \hat C_1
		 := 
		2\Big(1 \vee M \vee 2^{\varrho-1} 
			(M \sqrt{d})^\varrho 
			\Big(
				C_0 
				+ \|\mu \|^\varrho T^{\frac{\varrho}{2}} C_{BDG} C_0 \big(\bar \lambda_{(\sigma \sigma^{\top})^{-1}} \big)^{\frac{\varrho}{2}}
			\Big) \Big),
	$$
	and
	$$
		\hat C_{2,\varrho}\le \hat C_2 := C_{\rm bdg}~C_0 ~\big(\bar \lambda_{(\sigma \sigma^{\top})^{-1}} \big)^{\frac{\varrho}{2}}.
	$$
	\end{Lemma}
	\proof
		Let $\tilde b \in \R^d$ be a fixed vector. Set 
		  $Q^{t,x}  := \nabla X^{t,x} \tilde b$.	Then, it follows from direct computations that
	\b*
		\E \Big[ \max_{[t , t_{i+1}]} \|Q^{t,x}\|^\varrho \Big]\le
		C_0 \|\tilde b \|^\varrho.
	\e*
	
	Further, remember that  each $b_{q}$ is assumed to be bounded by $1$, so that 
	$\|b_q \sigma^{-1} \|^2$ is uniformly bounded by $\bar \lambda_{(\sigma \sigma^{\top})^{-1}}$. Then, direct  computations lead to
	$$
		C_{1,\varrho}
		 \le
		 1 \vee M^{\varrho} \vee 2^{\varrho-1} \Big(
				C_0 
				+ \|\mu \|^q T^{\frac{\varrho}{2}} C_{\rm bdg} C_0 \big(\bar \lambda_{(\sigma \sigma^{\top})^{-1}} \big)^{\frac{\varrho}{2}}
			\Big),
	$$
	and
	$$
		C_{2,\varrho}
		\le
		C_{\rm bdg} ~C_0 ~\big(\bar \lambda_{(\sigma \sigma^{\top})^{-1}} \big)^{\frac{\varrho}{2}}.
	$$

	It remains to use our specific choice of $\rho$ and $(p_{\ell})_{\ell \in L}$ in Proposition \ref{prop: representation Ynm} to conclude.
	\qed

	\vspace{0.5em}

	Let us now choose  $h'_{\circ}$ and $M_{h'_{\circ}}$ such that
	\be \label{eq:def_delta}
		h'_{\circ} ~<~ 1 \wedge \frac{ \hat C_1^{-(|L|-1)}}{(|L|+1)(|L| -1) \hat C_2},
	\ee
	and
	\begin{equation} \label{eq:def_M_delta}
		(M_{h'_\circ})^4
		~:=~ 
		\big(  \hat C_1^{1 - |L|} - h'_{\circ} ( |L| + 1) (|L| - 1) \hat C_2 \big)^{(1-|L|)^{-1}}.
	\end{equation}

\begin{Lemma} \label{lemm:estim_V_DV} Let the conditions of   Proposition \ref{prop: representation Ynm} hold. Then, 
	the ordinary differential equation $\eta'(t) = \sum_{\ell  \in L} \hat C_2 \eta(t)^{\|\ell\|_{1}}$
	with initial condition $\eta(0) = \hat C_1 \ge 1$
	has a unique solution on $[0, h'_{\circ}]$, and it is bounded by $(M_{h'_{\circ}})^2$.
	Moreover, 
	\be \label{eq:branch_bounded_by_Mh}
		\E\Big[
		\Big( \prod_{k \in \Kc_t} \frac{\hat C_1}{\bar F(t - T_{k-})} \Big)
		\Big( \prod_{k \in \bar \Kc_t \setminus \Kc_t} \frac{ \hat C_2}{p_{\xi_k} \rho(\delta_k)} \Big)
		\Big] 
		~\le~
		\eta(t) ~\le~ (M_{h'_\circ})^4,
	\ee
	for all $t \in [0, h'_{\circ}]$.
\end{Lemma}
\proof The result follows from exactly the same arguments as in    \cite[Lemma A.1]{BouchardTanZou}. 
\qed
\vs2

We can now conclude the proof of Proposition \ref{prop: representation Ynm}.
\vs2

{\bf Proof of Proposition \ref{prop: representation Ynm}.}
In view of \eqref{eq: majo puissance 4}, Lemma \ref{lemm:estim_V_DV} implies that $\{|U^{m}_{t,x}|^{2}+\|V^{m}_{t,x}\|^{2}, (t,x)\in [t_{i},t_{i+1})\x \R^{d}\}$ is uniformly integrable (with a bound that does not depend on $i<N_{h}$ for $0<h\le h'_{\circ}$ ). Then, arguing exactly  as in  \cite[Proposition A.2]{BouchardTanZou} leads to $\tilde u^{m}=\bar u^{m}$ on $[t_{i},t_{i+1})$. Combined with  \cite[Proposition 3.7]{HOTTW}, the uniform integrability also implies that  $D\tilde u^{m}=\tilde v^{m}\sigma$ on  $ (t_{i},t_{i+1})\x \R^{d}$, and one can  conclude 
from Theorem \ref{thm: main} that $\tilde v^{m}=\tilde u^{m}$ on  $(t_{i},t_{i+1})\x \R^{d}$. By the induction hypothesis of the beginning of this section, this proves that the statements of   Proposition \ref{prop: representation Ynm} hold. 
\qed

	\begin{Remark} 
		 The constants $\hat C_1$ and $\hat C_2$ (and hence $h'_{\circ}$ and $M_{h'_{\circ}}$) are clearly not  optimal for applications.
		For instance, if $\sigma \equiv \sigma_\circ$, for some non-degenerate constant matrix, the constants $C_{1,\varrho}$ and $C_{2,\varrho}$ can be significantly  simplified as shown in  \cite[Remark 3.9]{HOTTW}.
	\end{Remark}

\section*{Acknowledgements}

This work has benefited from the financial support of the Initiative de Recherche ``M\'ethodes non-lin\'eaires pour la gestion des risques financiers'' sponsored by AXA Research Fund.

\vspace{1mm}

Bruno Bouchard and Xavier Warin acknowledge the financial support of ANR  project CAESARS (ANR-15-CE05-0024)..

\vspace{1mm}

Xiaolu Tan  acknowledges the financial support of the ERC 321111 Rofirm, the ANR Isotace (ANR-12-MONU-0013), and the Chairs Financial Risks (Risk Foundation, sponsored by Soci\'et\'e G\'en\'erale) and Finance and Sustainable Development (IEF sponsored by EDF and CA).

\bibliographystyle{unsrt}

\end{document}